\documentclass[twocolumn]{autart}    

\usepackage{graphicx}          
\DeclareFontFamily{OT1}{pzc}{}
\DeclareFontShape{OT1}{pzc}{m}{it}{<-> s * [1.200] pzcmi7t}{}
\DeclareMathAlphabet{\mathpzc}{OT1}{pzc}{m}{it}
\usepackage{amsmath,amssymb}
\usepackage{graphicx}
\usepackage{algorithm}
\usepackage{algpseudocode}
\usepackage{fancybox}
\usepackage{bm}
\usepackage{amsfonts}
\usepackage{color}
\usepackage{textcomp}
\usepackage{subcaption}
\usepackage{caption}

\usepackage{mathrsfs}
\newcommand{\bfu}{\mathbf{u}}
\newcommand{\bfw}{\mathbf{w}}
\newcommand{\bfv}{\mathbf{v}}
\newcommand{\bfn}{\mathbf{n}} 
\newcommand{\bfe}{\mathbf{e}}
\newcommand{\M}{\mathcal{M}}
\renewcommand{\S}{\mathcal{S}}

\newcommand{\R}{\mathbb{R}}
\newcommand{\E}{\mathbb{E}}
\newcommand{\D}{\mathfrak{D}}
\newcommand{\calO}{\mathcal{O}}
\newcommand{\proj}{\mathpzc{Proj}}
\newcommand{\llangle}{\left\langle}
\newcommand{\rrangle}{\right\rangle}

\newtheorem{theorem}{Theorem}
\newtheorem{assumption}{Assumption}
\newtheorem{lemma}{Lemma}

\begin{document}

\begin{frontmatter}

\title{Geometry-preserving Numerical Scheme for Riemannian Stochastic Differential Equations} 

\thanks[footnoteinfo]{Victor Solo is the corresponding author.}

\author[unsw]{Xi Wang}\ead{xi.wang14@unsw.edu.au},    %
\author[unsw]{Victor Solo}\ead{v.solo@unsw.edu.au}, 
\address[unsw]{School of Electrical Engineering \& Telecommunications, UNSW, Sydney, Australia.}


\begin{abstract}

Stochastic differential equations (SDEs) on Riemannian manifolds have numerous applications in system identification and control. However, geometry-preserving numerical methods for simulating Riemannian SDEs remain relatively underdeveloped. In this paper, we propose the Exponential Euler–Maruyama (Exp-EM) scheme for approximating solutions of SDEs on Riemannian manifolds. The Exp-EM scheme is both geometry-preserving and computationally tractable. We establish a strong convergence rate of \(\mathcal{O}(\delta^{\frac{1 - \epsilon}{2}})\) for the Exp-EM scheme, which extends previous results obtained for specific manifolds to a more general setting. Numerical simulations are provided to illustrate our theoretical findings.

\end{abstract}

\end{frontmatter}

\section{Introduction}
Stochastic differential equations (SDEs), which model the state dynamics of systems influenced by random noise, play a vital role in a wide range of applications across fields such as engineering, physics, and finance~\cite{mao2007stochastic}. To analyze these systems in practice, numerical schemes that generate approximate solution through time discretization are essential. In Euclidean spaces, such methods have been extensively studied in the literature, including in~\cite{platen1999introduction} and~\cite{higham2001algorithmic}.

Recent research has focused on SDEs where the state variables are subject to geometric constraints. These constraints, which can be modeled by a Riemannian manifold, naturally arise in various applications. For example, in molecular dynamics, the Hamiltonian of the phase space must remain constant~\cite{tuckerman2000understanding}. Similarly, in aerospace control, the attitude of an aircraft must lie in the special orthogonal group~\cite{goran2016engineer}.

The growing number of applications involving Riemannian SDEs underscores the need to develop geometry-preserving numerical schemes (GPNS) for SDEs on general manifolds \( \M \). A GPNS ensures that the entire trajectory of the numerical solution of the SDE remains on the manifold. In contrast, most conventional Euclidean numerical schemes are not geometry-preserving, as they do not account for the underlying geometric constraints.

Several previous works~\cite{solo2021convergence,muniz2022higher,solo2024stratonovich} proposed GPNS for SDEs on specific manifolds, such as spheres, Stiefel manifolds, and matrix Lie groups. However, the algorithms and convergence analyses in these studies are heavily dependent on the particular structures of the manifolds. Consequently, these methods are not straightforwardly extendable to general Riemannian manifolds.

One methodology for constructing GPNS on general Riemannian SDEs is based on manifold projection~\cite{averina2019modification,ciccotti2008projection}, where each iterate is projected back onto the manifold \( \M \) after update. However, this approach faces several challenges. Due to the non-convexity of manifolds, the uniqueness of the projection is not guaranteed. Even when the projection map is well-defined, the computation of the projection typically involves solving an optimization problem without a closed-form solution~\cite{cherian2016Riemannian}, making the algorithm computationally intractable.

Recently, \cite{armstrong2022curved} proposed a numerical scheme for Riemannian SDEs based on the Castell–Gaines (CG) method and analysis the convergence order. The core idea of the geometric CG scheme is to locally approximate the SDE by a Euclidean ordinary differential equation (ODE). However, the resulting trajectory remains on the manifold only if the ODE is solved exactly. In practice, since the ODE is solved approximately, the resulting trajectories remain close to the manifold but do not strictly preserve the geometric structure.

Therefore, the problem of designing computationally tractable GPNS for SDEs on general manifolds remains open. Motivated by this research gap, we propose a GPNS for SDEs on general Riemannian manifolds based on the exponential map. The exponential map, often interpreted as ``moving along geodesics," offers a projection-free approach. Moreover, the exponential, which is directly operate on the manifold, ensures that the numerical solution preserves the geometric structure throughout the iteration.

The contributions of our paper are as follows:
\begin{itemize}
    \item We propose the Exponential Euler–Maruyama (Exp-EM) scheme for SDEs on general Riemannian manifolds. The Exp-EM scheme is geometry-preserving as well as computationally tractable\footnote{ As long as the exponential map is tractable, which is the case for many commonly used manifolds in practice, such as spheres, hyperbolic spaces, and matrix Lie groups~\cite{lee2018introduction}.}.
    \item We establish a strong convergence rate of \(\mathcal{O}(\delta^{\frac{1 - \epsilon}{2}})\) for the Exp-EM method, with arbitrarily small \(\epsilon > 0\). The convergence order matches the convergence rates achieved in previous works for specific manifolds and extends the results to general manifolds.
    \item We conduct numerical simulations to illustrate the effectiveness of our method.
\end{itemize}

The remainder of the paper is organized as follows. In Section II, we review the Riemannian geometry of embedded manifolds in Euclidean spaces. Section III presents the formal problem definition and discusses the convergence order for evaluating a GPNS. In Section IV, we propose the Exp-EM scheme for SDEs on a general manifold \( \M \). Section V provides the convergence analysis of the Exp-EM scheme. Section VI presents numerical simulations to demonstrate the effectiveness of our method, and Section VII concludes the paper. Detailed proofs are provided in the Appendix.

\section{Preliminaries}
In this section, we review embedded manifolds in Euclidean spaces and the Riemannian geometry of the manifold. For more details, we refer the reader to~\cite{lee2018introduction}.
\subsection{Manifold in Euclidean Space}

     A (embedded) manifold \( \M \subseteq \R^{n+m} \)  is defined as $$ \M := \left\{  x \in \R^{n+m} \mid h(x) = 0, h:\R^{n+m} \to \R^m \right\},$$ where $h$ is a smooth function with its Jacobi matrix $N_x := \left[ \nabla h^1(x), \nabla h^2(x), \dots, \nabla h^m(x) \right]$ has full rank \( m \). 
    
    Under this assumption, the columns of \( N_x \), i.e., the Euclidean gradients \( \{\nabla h^l(x)\},   l \in [m] \), span the normal space \( N_x\M \) with dimension $m$. Moreover, the tangent space \( T_x\M \) is defined as the orthogonal complement of \( N_x\M \) in \( \R^{m+n} \), i.e., $T_x\M = \left\{\bfv \in \R^{m+n} \mid N_x^T \bfv = 0 \right\}$.

    The projection of a vector \( \bfw \in \R^{n+m} \) onto the normal space \( N_x\M \) is given by the following linear map $$\bfw^{\perp} := N_x (N_x^T N_x)^{-1} N_x^T \bfw,$$ and projection of a vector \( \bfw \) onto \( T_x\M \) is given by $$\bfw^{\top}:= \bfw - N_x (N_x^T N_x)^{-1} N_x^T \bfw.$$
    \subsection{Riemannian Geometry of Embedded Manifold}
    Since \( \M \) is embedded within \( \R^{n+m} \), we naturally equip \( \M \) with the induced metric from the Euclidean inner product \( \llangle \cdot, \cdot \rrangle \). Specifically, we define a Riemannian metric \( g \) on \( \M \) by setting, for any \( x \in \M \) and \( \bfv_1, \bfv_2 \in T_x\M \), $g_x(\bfv_1, \bfv_2) = \llangle \bfv_1, \bfv_2 \rrangle$. This definition endows \( (\M, g) \) with a Riemannian structure, allowing us to study its geometric properties.

    The embedding geometry of \( \M \) in the Euclidean space is captured by the second fundamental form, which quantifies the extrinsic curvature of $\M$. For any tangent vector fields \( \bfu, \bfv \), the second fundamental form is given by
    \[
    \Pi(\bfu, \bfv) := (\nabla_\bfu \bfv)^\perp,
    \]
    where \( \nabla_\bfu \bfv :=\left( \sum_{i=1}^{m+n} u_i(x) \frac{\partial \bfv}{\partial x_i} \right) \) denotes the standard directional derivative in \( \R^{m+n} \). It can be shown that \( \Pi(\cdot,\cdot) \) is a symmetric bilinear map~\cite{lee2018introduction}.

A curve \( \gamma(t) \subset \M \) is a geodesic if its acceleration coincides with the second fundamental form, i.e, $$\nabla_{\dot\gamma} \dot \gamma = \Pi(\dot\gamma,\dot\gamma),$$ or equivalently, $\D(\dot\gamma,\dot\gamma) := (\nabla_{\dot\gamma} \dot \gamma )^\top = 0$. It is shown that the latter equation is uniquely determined by its initial position \( x \) and initial tangent vector \(\bfv \)~\cite{lee2018introduction}.

Accordingly, we define the exponential map \(\exp_x: T_x\M \to \M\), which maps a vector \( v \in T_x\M \) to the point $\exp_x(\bfv) = \gamma_{x,\bfv}(1),$ i.e., the point reached by following the geodesic starting at \(x\) with initial velocity \(\bfv\) for unit time. The exponential map plays a fundamental role in Riemannian geometry, establishing a connection between the tangent space and the global structure of \(\M\).

    \section{Problem Statement}
    In this section, we introduce the SDE evolving on Riemannian manifolds and formally state our problem.

    \subsection{SDE on Riemannian Manifolds}

    A Stratonovich SDE (S-SDE) evolving on a manifold \( \M = \{x \in \R^{n+m} \mid h(x) = 0\} \) can be written as~\cite{solo2024stratonovich}
    \[
    dx = \alpha_s(x) dt + \sum_{j=1}^d \alpha_j(x) \circ dW_j,
    \]
    where \( \alpha_s(x), \alpha_j(x) \in T_x\mathcal{M} \), \( W_j \) are independent standard Brownian motions with variance \( dt \), and \( \circ \) denotes the Stratonovich integral. Applying the Stratonovich-to-It\^o conversion~\cite{solo2024stratonovich}, we obtain the corresponding It\^o  SDE (I-SDE) as
\begin{align}\label{eq:manifold_sde}
\begin{cases}
    dx = \alpha_s(x) dt + \alpha_d(x) dt + \sum_{j=1}^d \alpha_j(x) dW_j, \\
    \alpha_d(x)  = \frac{1}{2} \sum_{j=1}^d \nabla_{\alpha_j} \alpha_j(x)
\end{cases}
\end{align}
The additional It\^o drifting term \( \alpha_d(x) \) can be decomposed as
\begin{align*}
    \alpha_d(x) = \frac{1}{2} \sum_{j=1}^d \D_{\alpha_j} \alpha_j + \frac{1}{2} \sum_{j=1}^d \Pi(\alpha_j, \alpha_j),
\end{align*}
where $ \D_{\alpha_j} \alpha_j  \in T_x\M$ and $ \Pi(\alpha_j, \alpha_j) \in N_x\M$.

The term \( \sum_{j=1}^d \Pi(\alpha_j, \alpha_j) \) reflects the influence of the curvature of manifold $\M$, leading to a normal drift in the SDE on \( \M \). The normal drift term, which tends to push the process away from the manifold, poses significant challenges in the development of geometry-preserving numerical schemes for SDEs on general manifold \( \M \).

\subsection{Problem Formulation}
In this paper, we aim to design geometry-preserving numerical schemes (GPNS) for SDEs on Riemannian manifold \( \M \). A GPNS discretizes a time interval \([0, T]\) into \( M \) time steps of length \( \delta = T / M \). At each time step, the scheme iteratively computes the state variable \( \hat{x}_t \) by incorporating the drift and diffusion terms, ensuring that the approximation \( \hat{x}(t) \) for \( t \in [0, T] \) remains on \( \M \).

The performance of a GPNS is characterized by its strong convergence order. A GPNS on \( \M \) is said to have a strong convergence order \( p \) if that exist a constant $C$ independent on $\delta$, such that the error satisfies  
\[
    \mathbb{E} \left[ \max_{0 \leq t \leq T} \| x(t) - \hat{x}(t) \|^2 \right]^{\frac{1}{2}} \leq C \delta^p.
\]  

\section{Exp-EM Scheme}
In this section, we introduce our Exponential Euler–Maruyama (Exp-EM) scheme. At each iteration \( m \), the Exp-EM scheme first applies a one-step Euler method to compute
\begin{align}\label{eq:def_wm}
    \bfw_m = \left( \alpha_s(\hat{x}_m) + \alpha_d(\hat{x}_m) \right) \delta + \sum_{j=1}^d \alpha_j \sqrt{\delta} \, \epsilon_{j,m}.
\end{align}
    
We then project \( \bfw_m \) onto the tangent space of the manifold, since the exponential map \(\exp_x: T_x \mathcal{M} \to \mathcal{M}\) only takes a tangent vector as input. The next iterate is then generated using the exponential map, i.e.,
\begin{align*}
    \begin{cases}
        \bfv_m = \bfw_m^\top \\
        \;\;\;\;\;= \left( \alpha_s(\hat{x}_m) + \alpha_d(\hat{x}_m)^\top \right) \delta + \sum_{j=1}^d \alpha_j \sqrt{\delta} \, \epsilon_{j,m}. \\
        \hat{x}_{m+1} = \exp_{\hat{x}_m}(\bfv_m).
    \end{cases}
\end{align*}

The Exp-EM scheme is detailed in Algorithm~\ref{alg:exp-eu}. We emphasize that projection \( \bfw_m^\top \) is different from the manifold projection used in~\cite{averina2019modification,ciccotti2008projection}. The manifold projection $\proj_\M(x)$ projects the point $x$ directly onto the manifold $\M$, defined as $\proj_\M(x) = \arg\min_{y \in \M} d(x, y)$, which is generally computationally intractable. In contrast, projection \( \bfw_m^\top \) is from the Euclidean space \( \mathbb{R}^{n+m} \) onto the linear tangent space \( T_{\hat{x}_m} \M \), which is a computationally tractable linear operator. 

To analyze the convergence of the Exp-EM scheme, we impose the following assumptions.

\begin{assumption}\label{asm:start}
    The coefficients $\alpha_s(x), \alpha_d(x), \alpha_j(x)$ are $L$-Lipschitz continuous on $\mathcal M$, i.e., there exists a constant $L > 0$ such that
        \[
        \max_{\mathcal I \in \{s,d,j\}} \| \alpha_{\mathcal I}(x) - \alpha_{\mathcal I}(x') \| \le L \| x - x' \|, \quad \forall x, x' \in \mathcal M.
        \] 
\end{assumption}
\begin{assumption}
    The coefficients $\alpha_s(x), \alpha_d(x), \alpha_j(x)$ are bounded on $\mathcal M$, i.e., there exists a constant $C > 0$ such that $\max_{\mathcal I \in \{s,d,j\}} \| \alpha_s(x) \| \le C, \quad \forall x \in \mathcal M.$
\end{assumption}

\begin{assumption}\label{asm:middle}
    The gradient of constraint function \( h \) is lower bounded by $L_1$, i,e., there exist a constant \( L_1 > 0 \) such that for any \( x \in \M \) and $l \in [m]$, $\|\nabla h^l(x)\| \ge L_1$.
\end{assumption}
Assumption~\ref{asm:middle} excludes cases where \(\|\nabla h\| \to 0\) on the manifold \(\mathcal{M}\), thereby avoiding pathological curvature scenarios~\cite{mantegazza2011lecture}. This assumption ensures that the constraint function \(h\) remains well-behaved, which is essential for both theoretical analysis and numerical stability.

Assumption~\ref{asm:middle} holds automatically in many practical settings. For example, on compact manifolds, it is guaranteed because the Jacobi \(N_x\) has full rank. Similarly, when the SDE evolves within a bounded domain, the assumption is also naturally satisfied. Therefore, Assumption~\ref{asm:middle} is practically reasonable.

\begin{assumption}\label{asm:end}
The constraint \( h \) is second-order and third-order Lipschitz, i.e., there exist constants \( L_2, L_3 > 0 \) such that for any \( \bfw, \bfw' \in \R^{m+n} \) and $l\in[m]$,
 \begin{align*}
 \begin{cases}
     \| \nabla^2 h^l(\bfw,\bfw')\| = \|\nabla_{\bfw}\nabla_{\bfw'}h^l(x)\| \le  L_2\|\bfw\|\|\bfw'\|,\\
     \| \nabla^3 h^l(\bfw)\| = \| \nabla_{\bfw}\nabla_{\bfw}\nabla_{\bfw}h^l(x)\| \le  L_3\|\bfw\|^3 .
\end{cases}
\end{align*}
\end{assumption}
\begin{algorithm}[t]
    \caption{Exp-EM Scheme}
    \begin{algorithmic}[1]
    \State \textbf{Input:} Time step size $\delta$, final time $T$, initial state $R_0 \in \M$.
    \State Divide the interval $[0, T]$ into $M$ subintervals with length $\delta$.
    \For{$m = 0, \dots, M$}
        \State Compute $\alpha_s(\hat x_m), \alpha_j(\hat x_m)$;
        \State Compute $\alpha_d(\hat x_m) =\frac{1}{2} \sum_{j=1}^d \nabla_{\alpha_j} \alpha_j(\hat x_m)$;
        \State Draw $\epsilon_{j,m} \sim \mathcal{N}(0, 1)$ and compute
        \begin{align*}
            &\bfv_m=  \Big( \alpha_s(\hat x_m)) + \alpha_d(\hat x_m)^\top\Big) \delta + \sum_{j=1}^d \alpha_j \sqrt{\delta} \epsilon_{j,m};
        \end{align*}
        \State Update $\hat x_{m+1} = \exp_{\hat x_m}(\bfv_m) $;
    \EndFor
    \State \textbf{Output:} Trajectory $\hat x(t) = \hat x_m, (m-1)\delta < t \le m\delta$.
    \end{algorithmic}
    \label{alg:exp-eu}
\end{algorithm}
By the above assumptions, we state the following lemma.
\begin{lemma}\label{lem:sec-fun}
    Under Assumptions~\ref{asm:middle} and~\ref{asm:end}, for any tangent vector fields \( \bfu, \bfv \), there holds
    \begin{align*}
        \|\Pi(\bfu,\bfv)\| \le C_2 \|\bfu\|\|\bfv\|,
    \end{align*}
    where \( C_2 := \sqrt{m} \frac{L_2}{L_1} \).
\end{lemma}
The proof of Lemma~\ref{lem:sec-fun} is provided in Appendix. Lemma~\ref{lem:sec-fun} provides a uniform bound on the second fundamental form, indicating that the normal drift in~\eqref{eq:manifold_sde} is controlled by the regularity of the manifold $\M$.

Using Lemma~\ref{lem:sec-fun}, Theorem~\ref{thm:exp} establishes a third-order approximation of the exponential map on the Riemannian manifold \( \M \).

\begin{theorem}\label{thm:exp}
    Under Assumptions~\ref{asm:middle} and~\ref{asm:end}, for any \( x \in \M \) and \( \bfv \in T_x\M \), the exponential map \( \exp_x(\bfv) \) satisfies
    \[
        \left\| \exp_x(\bfv) - x - \bfv - \frac{1}{2} \Pi(\bfv, \bfv) \right\| \leq C_3 \|\bfv\|^3,
    \]
    where \( C_3:=\sqrt{  C_2^4+ m (\frac{L_3+C_2L_2}{L_1})^2  } \).
\end{theorem}

The proof of Theorem~\ref{thm:exp} is provided in Appendix. Theorem~\ref{thm:exp} demonstrates that, the first-order term of the exponential map $\exp_x(\bfv)$ is \( \bfv \), which represents the primal tangential component of the exponential map. Moreover, the second-order term is \( \Pi(\bfv, \bfv) \), which corresponds to the primal normal component of the exponential map. The approximation plays a crucial role for the following convergence analysis.  

\section{Convergence Analysis}

With the help of the third-order approximation of the exponential map, we can carry out our convergence analysis. We first introduce the immediate sequence $ x_{\delta}$ as
\begin{align*}
    x_{\delta}(t) = x_0 +  \int_0^{t} (\alpha_s + \alpha_d)(\hat x_s) ds + \sum_{j=1}^d \int_0^{t}  \alpha_j(\hat x_s) dW_j.
\end{align*}
In this case, the error of our Exp-EM method, $$\mathbb{E} \left[ \| \hat{x}(t) - x_{\delta}(t) \|^2 \right],$$ can be controlled by the expectations
\begin{align*}
    a_\delta=\mathbb{E} \left[ \| x_{\delta}(t) - \hat{x}(t) \|^2 \right],\;\; b_\delta = \mathbb{E} \left[ \| \hat{x}(t) - x(t) \|^2, \right]
\end{align*}
which are analyzed in the following lemmas.
\begin{lemma}\label{lem:2}
    Under Assumptions~\ref{asm:start}-\ref{asm:end}, the error between $\hat x(t)$ and the $x_{\delta}(t)$ satisfies
    \begin{align*}
        \E \Big[ \sup_{0 \le t \le T}\| \hat x(t) - x_{\delta}(t) \|^2  \Big] \le \calO(\delta^{1-\epsilon}),
    \end{align*}
    for any arbitrary small $\epsilon > 0$.
\end{lemma}
The proof of Lemma~\ref{lem:2} is given in Appendix. 
\begin{lemma}
    Under Assumptions~\ref{asm:start}-\ref{asm:end}, there satisfires
    \begin{align*}
        \E  \Big[ \sup_{0 \le t \le T}\| x_{\delta}(t) - x(t) \|^2  \Big] \le \calO(\delta^{1-\epsilon}),
    \end{align*}
    for any arbitrary small $\epsilon > 0$.
\end{lemma}
\noindent \textit{Proof:} The proof is a straightforward generalization of the proof of Result IIIb in~\cite{solo2024stratonovich}. We omit the details here due to space limitations. \hfill $\square$ 

Combining the two lemmas, we derived our main result.
\begin{theorem}~\label{thm:cov-rate}
    Under Assumptions~\ref{asm:start}-\ref{asm:end}, the Exp-EM scheme converges to the real solution $x(t)$ of the SDE~\eqref{eq:manifold_sde} on $\M$ with strong order $\frac{1-\epsilon}{2}$, i.e.,
    \begin{align*}
        \E  \Big[ \sup_{0 \le t \le T}\| \hat x(t) - x(t) \|^2  \Big] \le \calO(\delta^{1-\epsilon}),
    \end{align*}
    for any arbitrary small $\epsilon > 0$.
\end{theorem}

Theorem~\ref{thm:cov-rate} establishes the convergence rate \(\mathcal{O}(\delta^{\frac{1 - \epsilon}{2}})\) of the Exp-EM scheme, which matches the rates achieved in previous works~\cite{solo2021convergence,solo2024stratonovich} for spheres and special orthogonal groups, and also being applicable to general manifolds.

\section{Numerical Simulation}

In this section, we present numerical simulations of Brownian motion on the sphere \(\S^n = \{ x \in \mathbb{R}^{n+1} \mid x^T x = 1 \}\) to demonstrate the effectiveness of our proposed Exp-EM scheme. We begin by introducing Brownian motion on \(\S^n\), and then conduct experiments in two cases: (i) \(n = 20\) to illustrate typical behavior in a moderate-dimensional setting, and (ii) \(n = 2000\) to evaluate performance in high-dimensional settings. We also compare our Exp-EM scheme with the classical Euclidean Euler–Maruyama (Eu-EM) scheme~\cite{higham2001algorithmic} and the geometric Castell–Gaines (G-CG) scheme~\cite{armstrong2022curved}.

\subsection{Brownian Motion on Spheres}
The Brownian motion on \( \S^n \) can be described by the following S-SDE
\[
    dx = \sum_{j=1}^{n+1} P(x)\bfe_j \circ dW_j,
\]
where \( P(x) = (I - xx^T) \) denotes the projection onto the tangent space \( T_x \S^n \), and \( \bfe_j \) is the \( j \)-th standard basis vector of \( \mathbb{R}^{n+1} \). The corresponding I-SDE of the Brownian motion is given by
\begin{align}\label{eq:BW-ISDE}
    dx = -\frac{n}{2} x\, dt + \sum_{j=1}^{n+1} P(x)\bfe_j \, dW_j,
\end{align}
where $\alpha_s =0$, $\alpha_d= -\frac{n}{2} x$ and $\alpha_j =  P(x)\bfe_j$. Applying It\^o’s formula, we obtain for any function \( f : \S^n \to \mathbb{R} \), 
\begin{align}\label{eq:def_ST}
    & f(x_t)  = S_t :=  f(x_0) + \sum_{j=1}^{n+1}\int_0^t   \llangle \nabla f(x), P(x)\bfe_j  \rrangle dW_{j} \nonumber \\ 
    & + \int_0^t \left[-\frac{n}{2}\llangle \nabla f(x ,x\rrangle + \frac{1}{2}{\rm tr}\left(\nabla^2 f(x) P(x)\right)\right] ds.
\end{align}

\subsection{Experiment in Moderate-dimensional Setting}
We begin our experiments by evaluating the performance of our Exp-EM scheme for solving~\eqref{eq:BW-ISDE} with \( n = 20 \), \( M = 1000 \), and varying the step size \( \delta \) from \( 2^{-10} \) to \( 2^{-3} \). The exponential map on the sphere \( \S^n \) is computed as $\exp_x(\mathbf{v}) = \cos(\|\mathbf{v}\|)\, x + \sin(\|\mathbf{v}\|)\, \frac{\mathbf{v}}{\|\mathbf{v}\|}$, showcasing the tractability of the Exp-EM scheme. For comparison, we also implement the Eu-EM scheme and the geometric G-CG scheme under the same conditions. As described in~\cite{armstrong2022curved}, the G-CG method employs a 4th-order Runge–Kutta (RK4) integrator.

\begin{figure*}[ht!]
\centering
\setkeys{Gin}{width=.95\linewidth}
    \begin{subfigure}[b]{0.32\linewidth}
        \centering
        \includegraphics{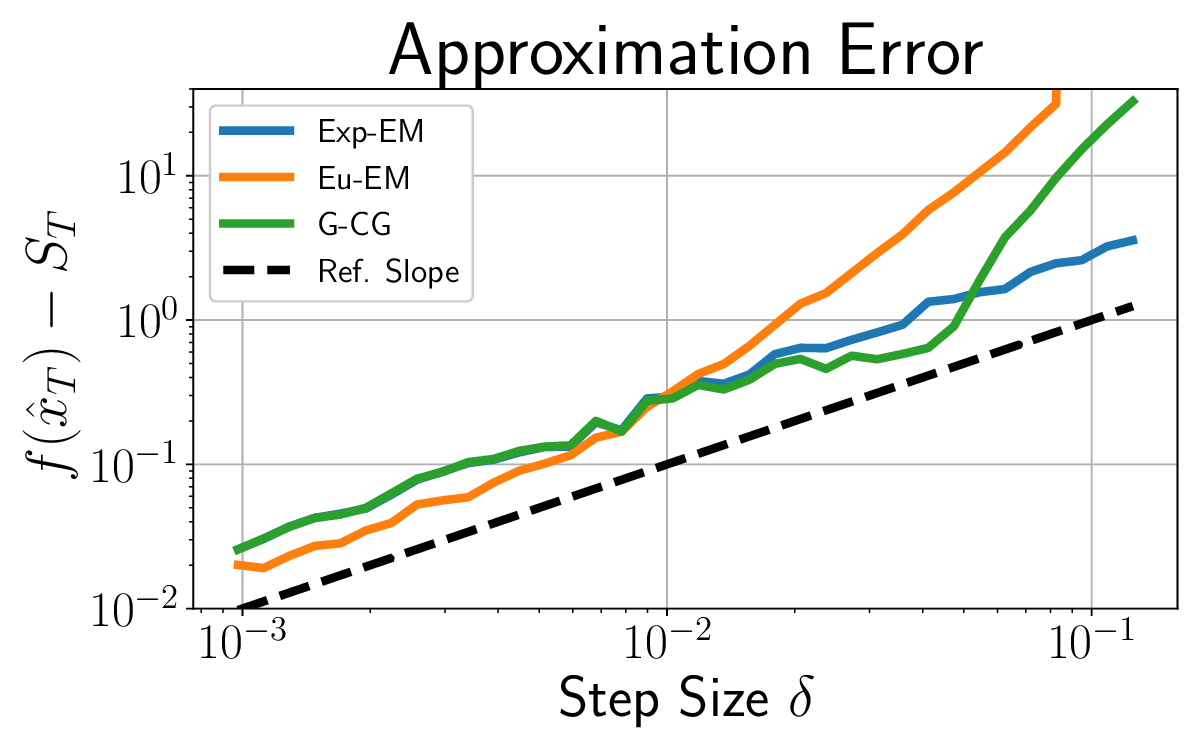}
        \caption{Approximation Error vs. Step Size}
        \label{fig:dim_20_err}
    \end{subfigure}
    \hfill
    \begin{subfigure}[b]{0.32\linewidth}
        \centering
        \includegraphics{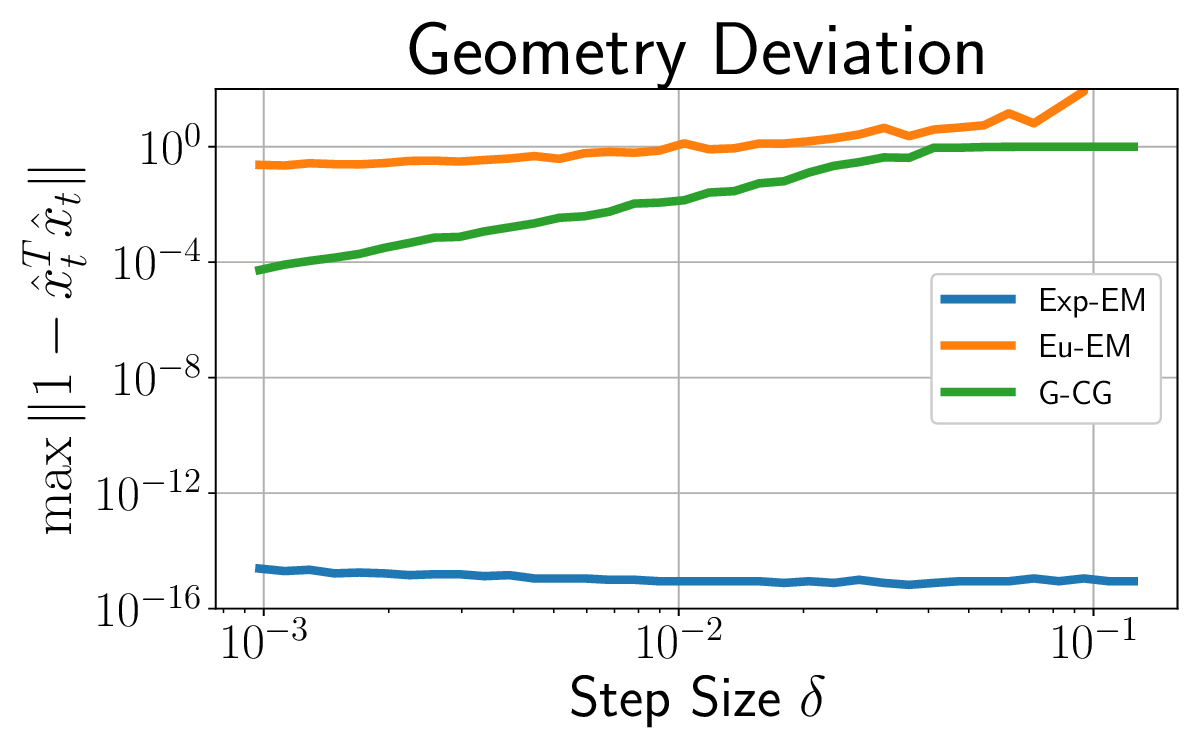}
        \caption{Geometry Deviation vs. Step Size} 
        \label{fig:dim_20_geo}
    \end{subfigure}
    \hfill
    \begin{subfigure}[b]{0.32\linewidth}
        \centering
        \includegraphics{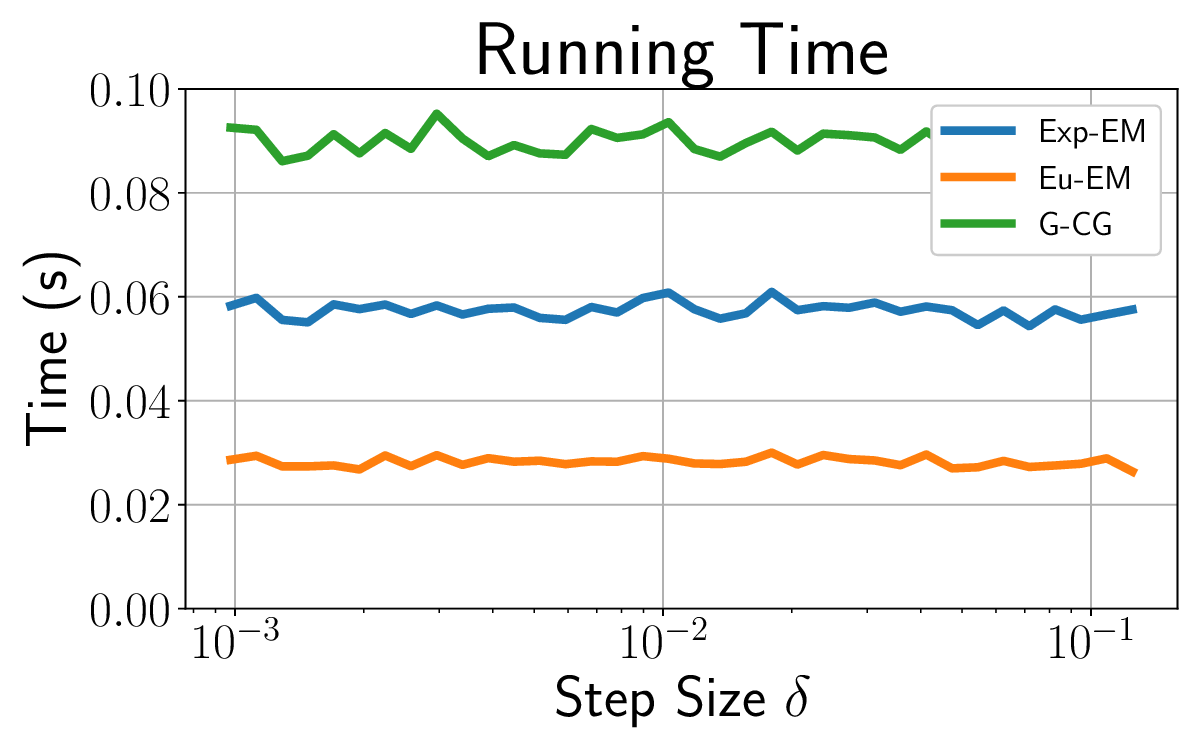}
        \caption{Running Time vs. Step Size}
        \label{fig:dim_20_time}
    \end{subfigure}
    \caption{Performance of Schemes for $n=20$}
    \label{fig:dim20}
\end{figure*}

We evaluate the performance of all the schemes in Fig.~\ref{fig:dim20}. In Fig.~\ref{fig:dim_20_err}, we plot the approximation error \( f(\hat x_T) - S_T \) versus the step size \( \delta \), where the function $f(x)= \frac{1}{2}x^2_{n+1}$ and the integral \( S_T \)~\eqref{eq:def_ST} is estimated using the Monte Carlo method by averaging over $50$ independent runs. The results in Fig.~\ref{fig:dim_20_err} demonstrate that our Exp-EM scheme exhibits an approximately \( \mathcal{O}(\delta) \) error across all tested step sizes. Moreover, both the Eu-EM and G-CG schemes also show an \( \calO(\delta) \) error when the step size \( \delta \) is small, but their approximation errors increase rapidly for larger step sizes. These results not only illustrate the correctness of our Exp-EM method, but also illustrate the robustness and stability of the Exp-EM method across a wider range of step sizes compared to the baseline.

We plot the geometric deviation \(\max\|1 - \hat{x}(t)^T \hat{x}(t)\|\) versus the step size \(\delta\) in Fig.~\ref{fig:dim_20_geo}. The Eu-EM method fails to preserve the geometric constraint with a deviation of \(\calO(10^0)\) for all step sizes. The G-CG method stays closer to the manifold, showing a deviation of \(\calO(10^{-4})\) for small \(\delta\), and the deviation increases to \(\mathcal{O}(10^0)\) as \(\delta\) approches \(\mathcal{O}(10^{-1})\). In contrast, our Exp-EM method maintains a geometric deviation of \(\mathcal{O}(10^{-15})\) across all step sizes, demonstrating the geometry-preserving property.

We also report the running times in Fig.~\ref{fig:dim_20_time}. The Eu-EM method is the fastest among all algorithms. Our Exp-EM method runs faster than the G-CG scheme, as the RK4 in the G-CG requires applying the projection \(P(x)\) four times per iteration. These above results demonstrate that our Exp-EM method reaches an $\calO(\delta)$ error and achieves geometry preservation with relatively low computational cost.

\begin{figure*}[ht!]
\centering
\setkeys{Gin}{width=.95\linewidth}
    \begin{subfigure}[b]{0.32\linewidth}
        \centering
        \includegraphics{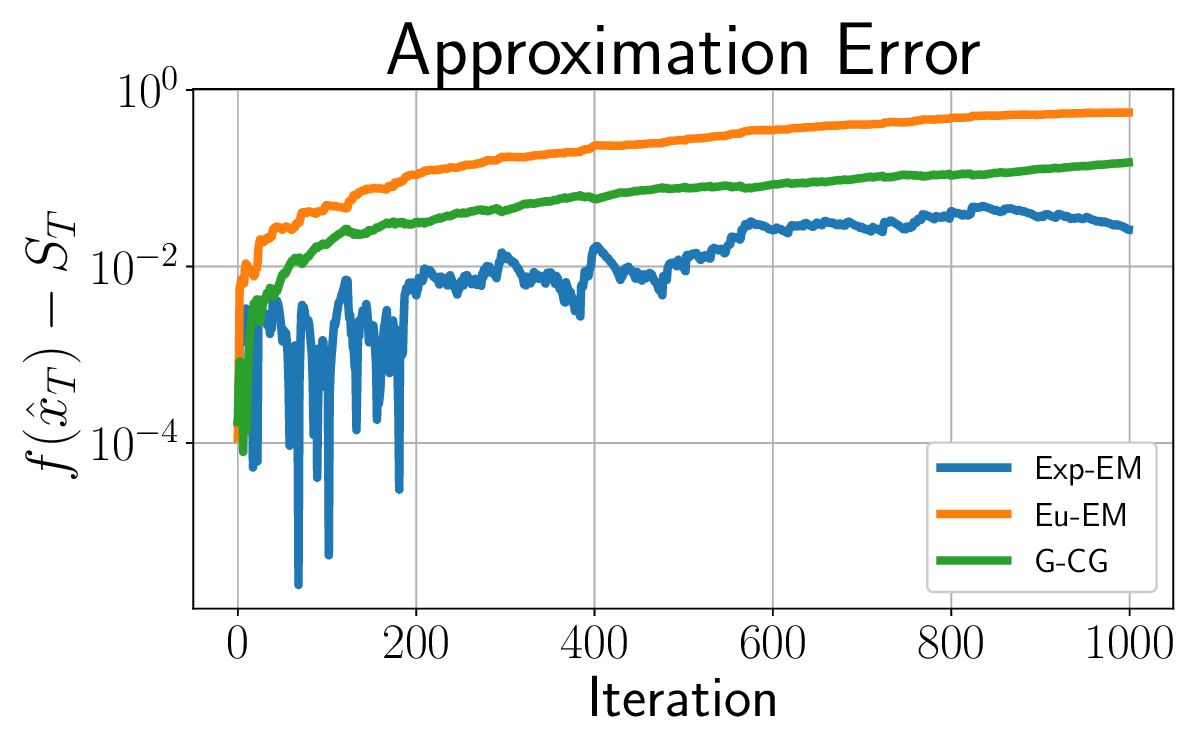}
        \caption{Approximation Error vs. Iteration}
        \label{fig:dim_2000_err3}
    \end{subfigure}
    \hfill
    \begin{subfigure}[b]{0.32\linewidth}
        \centering
        \includegraphics{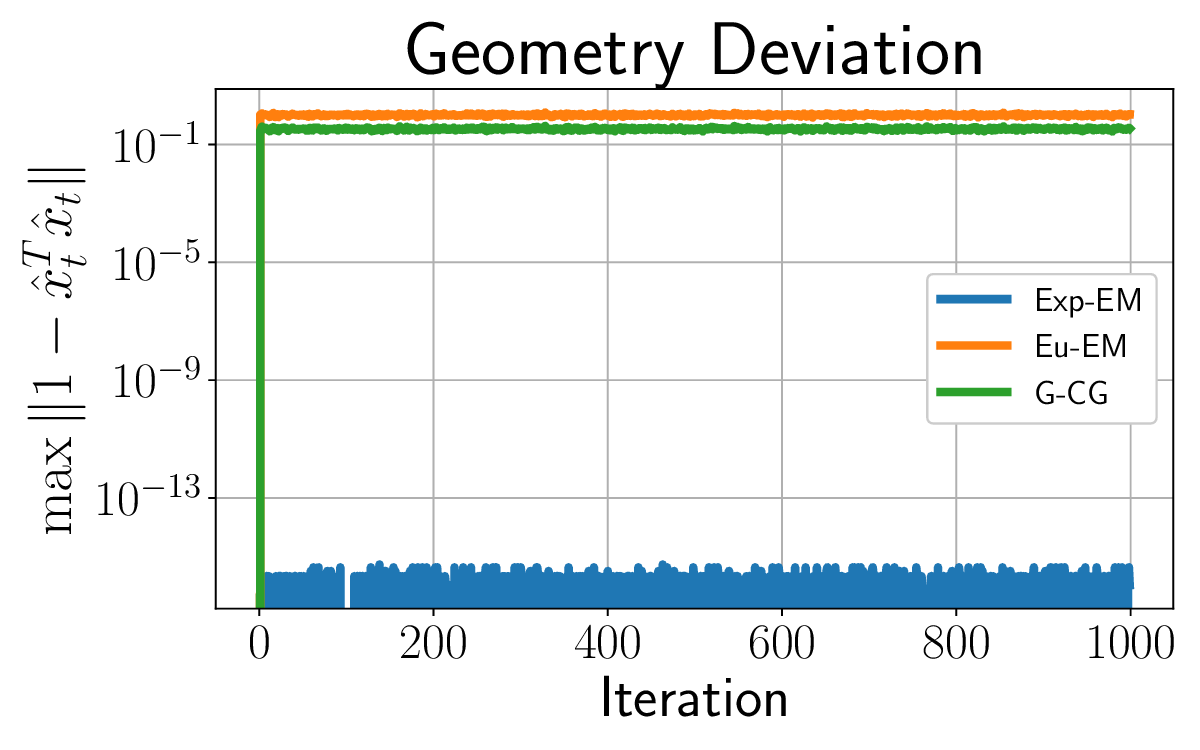}
        \caption{Geometry Deviation vs. Iteration} 
        \label{fig:dim_2000_geo3}
    \end{subfigure}
    \hfill
    \begin{subfigure}[b]{0.32\linewidth}
        \centering
        \includegraphics{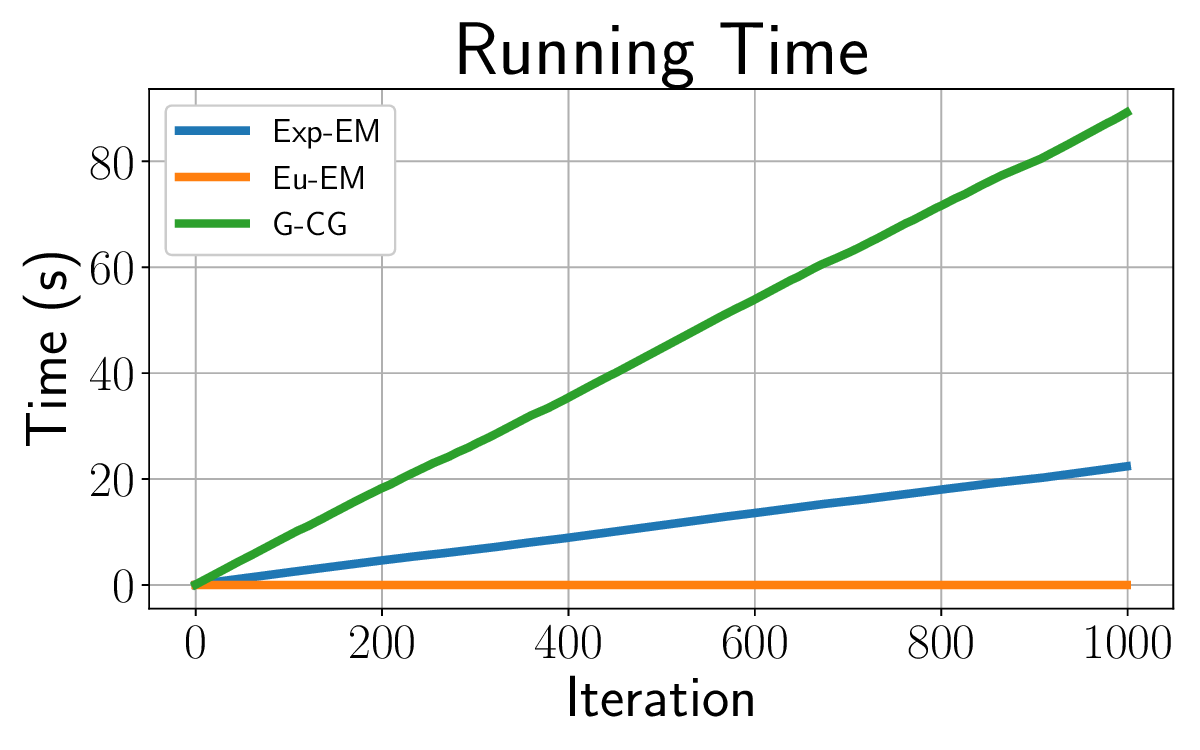}
        \caption{Running Time vs. Iteration}
        \label{fig:dim_2000_time3}
    \end{subfigure}
    \caption{Performance of Schemes for $n=2000$ and $\delta = 0.001$}
    \label{fig:dim2000_1e-3}
\end{figure*}

\begin{figure*}[ht!]
\centering
\setkeys{Gin}{width=.95\linewidth}
    \begin{subfigure}[b]{0.32\linewidth}
        \centering
        \includegraphics{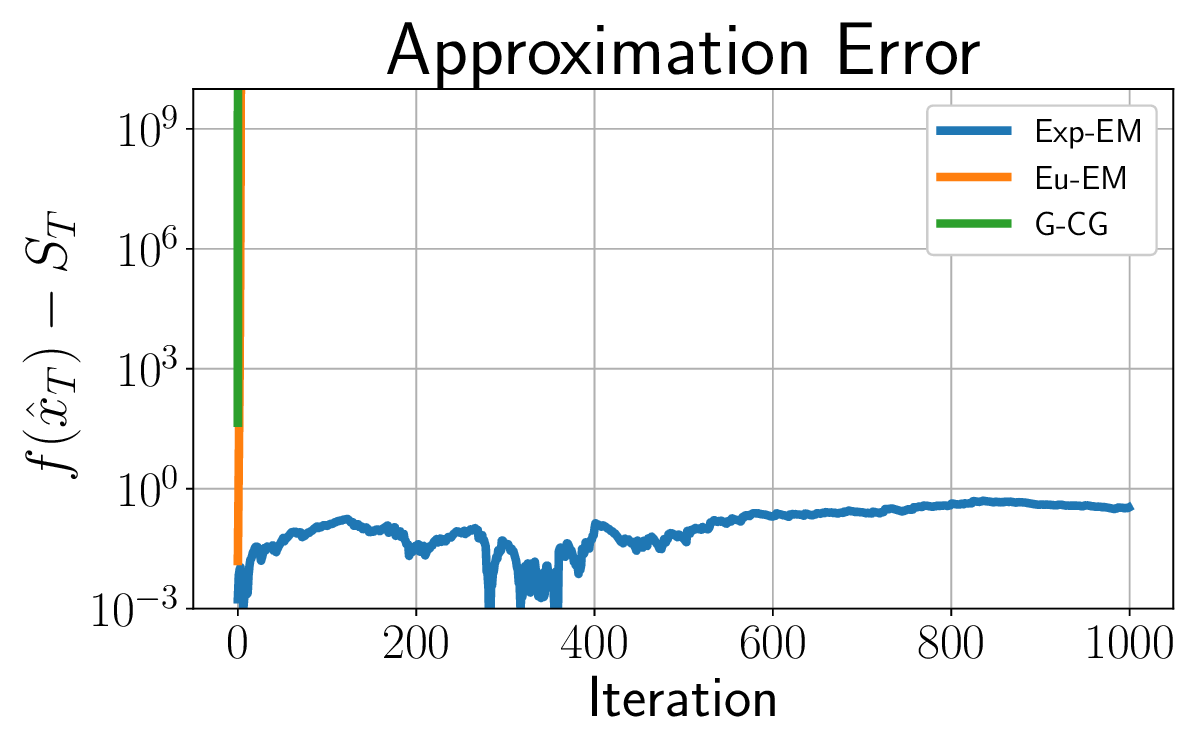}
        \caption{Approximation Error vs. Iteration}
        \label{fig:dim_2000_err2}
    \end{subfigure}
    \hfill
    \begin{subfigure}[b]{0.32\linewidth}
        \centering
        \includegraphics{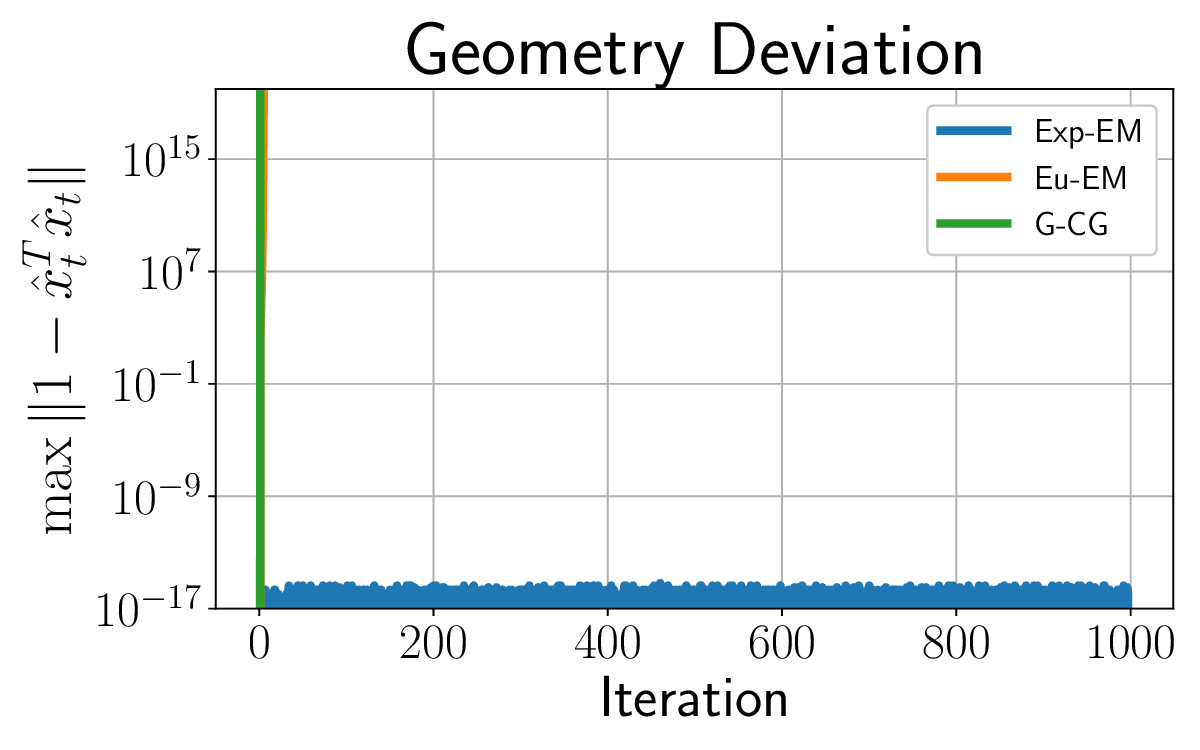}
        \caption{Geometry Deviation vs. Iteration} 
        \label{fig:dim_2000_geo2}
    \end{subfigure}
    \hfill
    \begin{subfigure}[b]{0.32\linewidth}
        \centering
        \includegraphics{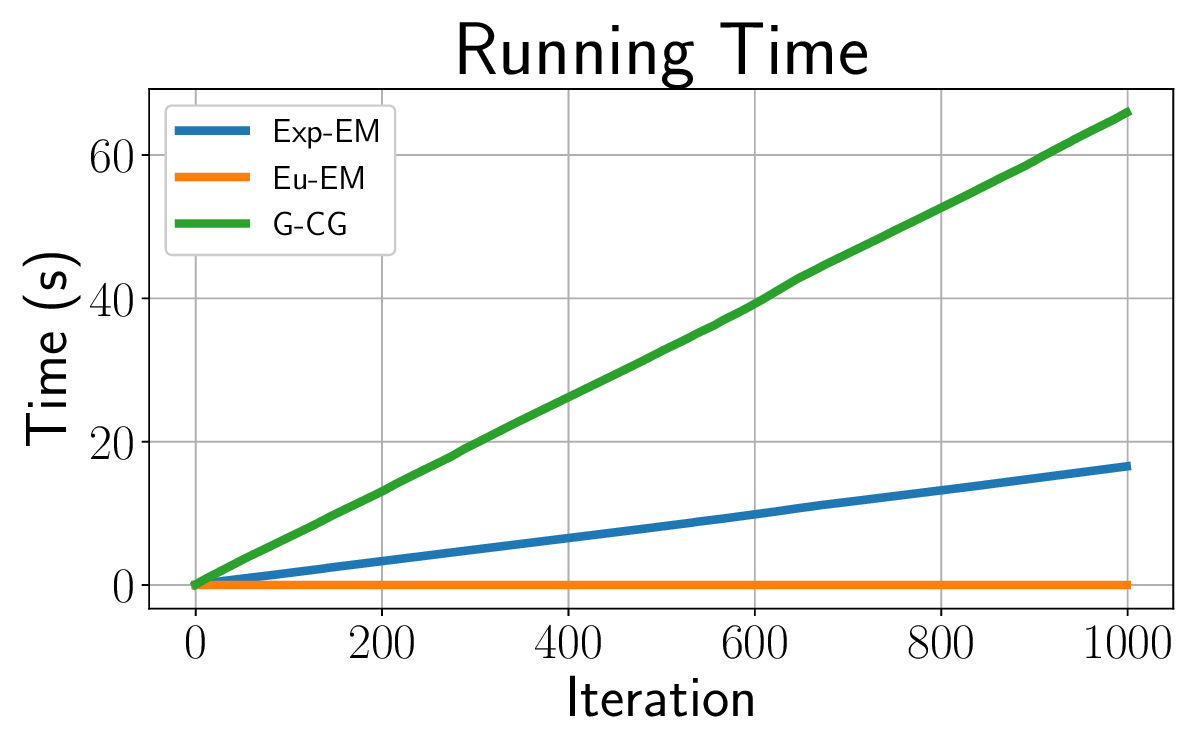}
        \caption{Running Time vs. Iteration}
        \label{fig:dim_2000_time2}
    \end{subfigure}
    \caption{Performance of Schemes for $n=2000$ and $\delta = 0.01$}
    \label{fig:dim2000_1e-2}
\end{figure*}

\subsection{Experiment in High-dimensional Setting}

In addition, we test all the above schemes for solving~\eqref{eq:BW-ISDE} with \( n = 2000 \) and \( M = 1000 \) to evaluate their effectiveness in high-dimensional settings. The performance of the schemes is presented in Fig.~\ref{fig:dim2000_1e-3} and Fig.~\ref{fig:dim2000_1e-2} for step sizes \(\delta = 0.001\) and \(\delta = 0.01\), respectively.

As shown in Fig.~\ref{fig:dim2000_1e-3}, our Exp-EM method achieves an error of \(\mathcal{O}(10^{-2})\) with a geometric deviation of \(\mathcal{O}(10^{-17})\) when the step size $\delta=0.001$. Both the G-CG and Eu-EM methods yield errors around \(\mathcal{O}(10^{-1})\) and geometric deviations on the order of \(\mathcal{O}(10^{-1})\). These results indicate that the Exp-EM method provides more accurate solution while maintaining a relatively low computational cost, clearly demonstrating the advantage of our Exp-EM methods. Furthermore, in Fig.~\ref{fig:dim2000_1e-2}, both the Eu-EM and G-CG methods encounter numerical overflow when $\delta=0.01$, whereas our Exp-EM method continues to preserve the geometric constraint and maintains an error of \(\mathcal{O}(10^0)\). These results highlight the potential stability and robustness of the Exp-EM method in large-scale computations.

\section{Conclusion}
In this paper, we studied geometry-preserving numerical schemes (GPNS) for stochastic differential equations (SDE) on general manifolds. We proposed the Exponential Euler–Maruyama (Exp-EM) scheme, which preserves the geometric constraint of the manifold and is computationally tractable. We also established a strong convergence rate of \(\mathcal{O}(\delta^{\frac{1 - \epsilon}{2}})\) for the Exp-EM scheme. Numerical simulations illustrate our theoretical findings and demonstrate the advantages of our method compared to existing baselines.

In future work, we aim to improve the strong convergence rate to \(\mathcal{O}(\delta^{\frac{1}{2}})\) and to develop higher-order GPNS for SDEs over general manifolds.

\bibliographystyle{plain}        
\bibliography{autosam}           



\appendix
\section{Proof of Lemma~\ref{lem:sec-fun}}
    Since the second fundamental form \( \Pi(\bfu, \bfv) \in N_x\M \), it is sufficient to perform computations in the normal space. Specifically, we consider the inner product $\llangle \Pi(\bfu,\bfv) ,\nabla h^l \rrangle,$ for $l \in [m]$. Since \( \bfu \in T_x\M \), we have $\llangle \bfu, \nabla h^l \rrangle = 0.$ Differentiating both sides along the vector \( \bfv \), we obtain  
    \begin{align*}
       0 = & \llangle \nabla_\bfu \bfv, \nabla h^l \rrangle + \llangle  \bfv, \nabla_\bfu  \nabla h^l \rrangle\\
         = & \llangle \Pi(\bfu,\bfv), \nabla h^l \rrangle + \nabla^2 h^l(\bfu,\bfv).
    \end{align*}
    This implies that for unit vector $ \bfn^l = \frac{\nabla h^l}{\|\nabla h^l\|}$, we have  
    \begin{align*}
         \llangle \Pi(\bfu,\bfv), \bfn^l \rrangle \le \frac{\|\nabla^2 h^l(\bfu,\bfv)\|}{\|\nabla h^l\|} \leq \frac{L_2}{L_1} \|\bfu\| \|\bfv\|.
    \end{align*}
    Since the set  
    $ \left\{ \bfn^l \right\}, l\in[m] $  
    forms a unit basis for the normal space \( N_x\M \), we have
    \[
        \| \Pi(\bfu, \bfv) \| \leq \sqrt{m} \frac{L_2}{L_1} \|\bfu\| \|\bfv\|,
    \]
    which completes our proof. \hfill $\square$
\section{Proof of Theorem~\ref{thm:exp}}
We first introduce a lemma about the Lipschitz continuity of the second fundamental form.

\begin{lemma}\label{lemma:sec-lip}
    Under Assumptions 3–4, for any tangent vector fields $\bfv$, there holds
    \begin{align*}
      \left\|  \nabla_\bfv\Pi(\bfv,\bfv) \right\| \le C_3 \|\bfv\|^3,
    \end{align*}
    where $C_3:=\sqrt{ C_2^4+  m(\frac{L_3+C_2L_2}{L_1})^2 }$.
\end{lemma}

\noindent {\it Proof.}
For simplicity, we denote $\Pi_\bfv$ as $\Pi(\bfv,\bfv)$. We compute the tangent and normal components of \( \nabla_\bfv \Pi_\bfv \), respectively. For any tangent vector \( \bfu \in T_x\M \), we have $\llangle \Pi_\bfv , \bfu \rrangle = 0$. Differentiating both side along \( \bfv \) yields $ 0 =\llangle \nabla_\bfv\Pi_\bfv , \bfu \rrangle + \llangle \Pi_\bfv , \nabla_\bfv\bfu \rrangle$.
Rearranging the terms, we obtain  
\begin{align*}
    \left|\llangle \nabla_\bfv\Pi_\bfv , \bfu \rrangle\right| = \left| \llangle\Pi_\bfv, \Pi(\bfv, \bfu) \rrangle \right| \le C^2_2\|\bfv\|^3\|\bfu\|,
\end{align*}
which indicates that the tangent component satisfies  
\begin{align}\label{eq:lem4-tg}
    \|(\nabla_\bfv\Pi_\bfv)^\top\| \le  C^2_2\|\bfv\|^3.
\end{align}
Next, we compute the normal component \( (\nabla_\bfv\Pi_\bfv)^\perp \), for the unit normal basis \( \bfn^l = \frac{\nabla h^l}{\|\nabla h^l\|} \), we have
\begin{align*}
     \llangle \nabla_\bfv \Pi_\bfv, \bfn^l \rrangle =  & \frac{1}{\nabla h^l} \left(\nabla_\bfv\llangle \Pi_\bfv,\nabla h^l \rrangle -\llangle \Pi_\bfv,\nabla_\bfv\nabla h^l \rrangle  \right) \\
     = &\frac{1}{\nabla h^l}\left(\nabla^3h^l(\bfv)+ \nabla^2h(\Pi(\bfv,\bfv),-\bfv)\right) \\ 
      \leq & (\frac{L_3+C_2L_2}{L_1})  \|\bfv\|^3.
\end{align*}
Following the proof in Lemma~\ref{lem:sec-fun}, we obtain  
\begin{align}\label{eq:lem4-nor}
    \|(\nabla_\bfv\Pi_\bfv)^\perp\| \le  \sqrt{m} (\frac{L_3+C_2L_2}{L_1})\|\bfv\|^3.
\end{align}
Combining~\eqref{eq:lem4-tg} and~\eqref{eq:lem4-nor}, we conclude that  
\begin{align*}
    \|\nabla_\bfv\Pi(\bfv,\bfv)\| \le \sqrt{ C_2^4+ m(\frac{L_3+C_2L_2}{L_1})^2  } \|\bfv\|^3,
\end{align*}
which completes our proof. \hfill $\square$

\noindent \textit{Proof of Theorem~\ref{thm:exp}.}
 Denote $y = \exp_x(\bfv)$ and $\gamma(t) =\exp_x(t\bfv)$ as the geodesic connecting $x$ and $y$. We denote the tangent vector of $\gamma(t)$ as $\bfu(t)$. From the classical calculus in Euclidean spaces, we have $ y-x= \int_0^1 \bfu(t) dt$, and it follows that
    \begin{align*}
        y-x-\bfv = & \int_0^1 \bfu(t)- \bfu(0) dt = \int_0^1 \int_0^t  \frac{d}{ds}\bfu(s) ds dt.
    \end{align*}
    Since \( \bfu(s) \) lies along the geodesic \( \gamma(t) \), the derivative \( \frac{d}{ds} \bfu(s) \) is equal to the derivative \( \nabla_{\dot{\gamma}} \bfu(s) = \nabla_{\bfu} \bfu \). Because \( \gamma \) is a geodesic, the relationship $\nabla_{\bfu} \bfu  = \Pi(\bfu,\bfu)$ gives us
\begin{align*}
    \|& y - x - \bfv - \frac{1}{2} \Pi(\bfv, \bfv) \nonumber\| \\
    = & \int_0^1 \int_0^t \Pi(\bfu(s), \bfu(s)) - \Pi(\bfu(0), \bfu(0)) \, ds \, dt \nonumber \\
    = & \int_0^1 \int_0^t \int_0^s \nabla_\bfu \big( \Pi(\bfu(r), \bfu(r)) \big) \, dr \, ds \, dt \\
    \le & \frac{C_3}{6} \max_{r \in [0,1]} \|\bfu(r)\|^3 = \frac{C_3}{6} \|\bfv\|^3,
\end{align*}
The last two inequalities come from the Lipschitz continuity from Lemma~\ref{lemma:sec-lip} and the fact that the norm of \( \bfu \) remains constant along the geodesic. This is due to the face\( \nabla_\bfu \llangle \bfu, \bfu \rrangle = 2\llangle \Pi(\bfu, \bfu), \bfu \rrangle = 0 \). Here, we complete the proof. \hfill \( \square \)

\section{Proof of Lemma~\ref{lem:2}}
From the update rule, we have
    \begin{align*}
    \begin{cases}
        \hat x_m =  \exp_{\hat x_{m-1}}(\bfv_{m-1}) \\
                 \;\;\;\;\;\;=  \hat x_{m-1} + \bfv_{m-1} + \frac{1}{2} \Pi(\bfv_{m-1},\bfv_{m-1})+ E^a_{m-1} \nonumber \\
                 \;\;\;\;\;\;=  x_0 + \sum_{k=0}^{m-1} \left( \bfv_{k} + \frac{1}{2} \Pi(\bfv_{k},\bfv_{k}) + E^a_k \right), \\
             E^a_k = \exp_{\hat x_k}(\bfv_k) -\bfv_{k} - \frac{1}{2} \Pi(\bfv_{k},\bfv_{k}).
    \end{cases}
    \end{align*}
    From the definition of $\bfw_k$ in~\eqref{eq:def_wm}, we note that $\bfw_k - \bfv_k = \alpha_d(\hat x_k)^\perp= \frac{1}{2}\sum_{j=1}^d\Pi(\alpha_j,\alpha_j)$, and thus we have
    \begin{align*}
             \hat x_m - x_{\delta}(m\delta) = &\sum_{k=0}^{m-1} \Big(  \frac{1}{2} \Pi(\bfv_{k},\bfv_{k}) - \sum_{j=1}^d\Pi(\alpha_j,\alpha_j)+E^a_k \Big)  \\
          = & \sum_{k=0}^{m-1} \Big(E^a_{k}+E^{b}_{k} +  E^{c}_{k} +  E^{d}_{k}\Big).
    \end{align*}
    where
    \begin{align*}
    \begin{cases}
        E^{b}_{k} = \frac{1}{2} \delta   \sum_{j = 1}^d \Pi(\alpha_j,\alpha_j) (\epsilon_{j,k}^2 - 1) \\ 
        E^{c}_{k} = \frac{1}{2} \delta  \sum_{j' \neq j} \Pi(\alpha_j,\alpha_k) \epsilon_{j,k} \epsilon_{j',k}\\
        E^{d}_{k} = \frac{1}{2} \delta^{\frac{3}{2}} \sum_{j=1}^d \Pi(\alpha_s+\alpha_d^\top,\alpha_j)\epsilon_{j,k} \\ 
        \;\;\;\;\;\;\;\;\;\;\;\; + \frac{1}{2} \sum_{j=1}^d  \delta^2 \Pi(\alpha_s+\alpha_d^\top,\alpha_s+\alpha_d^\top)\epsilon_{j,k} 
    \end{cases}
    \end{align*}
    Then we have, for $t \in [m\delta, (m+1)\delta]$,
    \begin{align*}
        & \hat x(t) - x_{\delta}(t) = \hat x_m - x_{\delta}(m\delta) +  x_{\delta}(m\delta) - x_{\delta}(t) \\
        = &\sum_{k=0}^{m-1} \Big(E^a_{k}+E^{b}_{k} +  E^{c}_{k} +  E^{d}_{k}\Big) + (t - m\delta)  (\alpha_s+\alpha_d)(\hat x_{m}) \\ 
        & + \sum_{j=1}^d   \alpha_j(\hat x_{m}) \big( W_j(t)- W_j(m\delta) \big).
    \end{align*}
    Taking the expectation on both sides, we have
    \begin{align*}
        & \E \Big[ \sup_{0 \le t \le T}\| \hat x(t) - x_{\delta}(t) \|^2 \Big] \le 6\sum_{\mathcal I =a}^f \Delta_{\mathcal I},    
    \end{align*}
    where
    \begin{align*}
        \begin{cases}
        \Delta_{\mathcal{I}} = \E \left[ \sup_{m \in [M]} \|  \sum_{k=0}^{m-1} E^{\mathcal{I}}_{k} \|^2 \right], \mathcal I \in \{a,b,c,d\}, \\
        \Delta_e = \E \left[ \sup_{\mathcal T} \| (t - m\delta) (\alpha_s+\alpha_d)(\hat x_{m}) \|^2 \right], \\
        \Delta_f = \E \left[ \sup_{\mathcal T} C^2 \| \sum_{j=1}^d \big( W_j(t)- W_j(m\delta) \big) \|^2 \right]. \\
        \mathcal T = \{ m \in [M], t \in [m\delta, (m+1)\delta] \}. 
        \end{cases}
    \end{align*}
    From Lemma 1 in~\cite{solo2021convergence}, we have $\Delta_f=\calO(\delta^{1-\epsilon})$. Applying the boundness of the second fundamental form, we also obtain $\Delta_d=o(\delta)$ and $\Delta_e=\calO(\delta)$. Therefore, it remains to show that $\Delta_a,\Delta_b,\Delta_c=\calO(\delta)$.

    We begin by analyzing \(\Delta_b\). Since \( \epsilon_{j,k} \) are independent standard normal random variables, it follows that \(\mathbb{E}[(\epsilon_{j,k}^2 - 1)(\epsilon_{j',k}^2 - 1)] = 0\) for \( j \ne j' \), and \(\mathbb{E}[(\epsilon_{j,k}^2 - 1)^2] = 2\).  
    Therefore, we obtain
    \begin{align*}
        \mathbb{E} [ \|E_k^b\|^2 ] = & \frac{1}{4} \delta^2 \sum_{j=1}^{d} \|\Pi(\alpha_j,\alpha_j)\|^2 \mathbb{E} [ (\epsilon_{j,k}^2 -1)^2 ] \\
        \le & \frac{1}{4} \delta^2  2 \sum_{j=1}^{d} \|\Pi(\alpha_j,\alpha_j)\|^2 \leq \frac{1}{2} \delta^2  d (C_2 C^2)^2.
    \end{align*}
    Observe that the sum $\sum_{k=1}^{m} E_k^b$ forms an martingale. By Doob's maximal inequality, we have
    \begin{align}\label{eq:c1}
        \Delta_b= & \mathbb{E} \Big[ \sup_{ m \in [M]} \Big\|\sum_{k=0}^{m-1}  E^{b}_{k} \Big\|^2 \Big] \leq  4 \mathbb{E} [\sum_{k=0}^{M-1}\| E^{b}_{k} \|^2] \nonumber \\ 
         \leq & 2 \delta^2 d (C_2 C^2)^2 M = \delta d (C_2 C^2)^2 T = \calO(\delta).
    \end{align}
    Now we turn to the term \( \Delta_c \). Since \( \mathbb{E}[\epsilon_{j,k} \epsilon_{j',k}] = 0 \) for \( j \neq j' \), we have
    \begin{align*}
    \E\left[\|E_k^c\|^2\right] \le 3 \delta^2 \sum_{j \neq j'} \|\Pi(\alpha_j, \alpha_{j'})\|^2 \mathbb{E}[\epsilon_{j,k}^2 \epsilon_{j',k}^2].
    \end{align*}
    Note that the sum \( \sum_{k=1}^{m} E_k^c \) also forms a martingale. Applying Doob’s maximal inequality and the fact $\mathbb{E}[\epsilon_{j,k}^2 \epsilon_{j',k}^2]=1$, we get
    \begin{align}\label{eq:c2}
    \Delta_c \leq 4 \mathbb{E} \big[ \sum_{k=0}^{M-1} \|E_k^\delta\|^2 \big] \leq 12 d C_2 C^2 \delta^2 M  = \mathcal{O}(\delta).
    \end{align}
    The remaining part is to analyze $\Delta_a$. First, we have
    \begin{align*}
        \Delta_a =\E \Big[ \sup_{m \in [M]} \big\| \sum_{k=0}^{m-1}  E_k \big\|^2 \Big] 
        \le  M^2 \sup_{m\in [M]} \E \Big[  \| E_k \|^2 \Big].
    \end{align*}
    From Theorem~\ref{thm:exp} ,we have
    \begin{align*}
        \Delta_a \le &  M^2 \sup_{m\in [M]}  \E[ \| E_k \|]^2  \le C_3 M^2 \sup_{m\in [M]} \E[\| \bfv_k \|^6]\\
        \le &  C_3M^2 \big( C \delta + \sqrt{\delta} C \sum_{j=1}^d E[| \epsilon_{j,k} |] \big)^6 =M^2\calO(\delta^3),
    \end{align*}
    which implies that 
    \begin{align}\label{eq:c3}
        \Delta_a=M^2\calO(\delta^3) = T^2\calO(\delta).
    \end{align}
    Combining~\eqref{eq:c1}-\eqref{eq:c3}, we complete our proof. \hfill $\square$.
\end{document}